\documentclass{article}
\usepackage{graphicx}
\usepackage{pgfplots}
\usepackage{subfigure}
\usepackage{caption}
\usepackage[justification=centering]{caption}
\usepackage{amsmath}
\usepackage{amssymb}
\usepackage{amsthm}
\usepackage{listings}
\usepackage{epstopdf}
\usetikzlibrary{arrows}
\usepackage{tikz}
\usepackage{float}
\usepackage[linesnumbered,ruled]{algorithm2e}
\usetikzlibrary{shapes.geometric, arrows}
\usepackage[utf8]{inputenc}
\usepackage{enumitem}
\usepackage[top=2cm, bottom=2cm, left=2cm, right=2cm]{geometry}
\usepackage{multicol}

\newtheorem{theorem}{Theorem}

\usepackage[titletoc,title]{appendix}
\usepackage{array}
\usepackage{authblk}

\usepackage[compress]{cite}

\begin{document}

	\title{Abstract Fractals}	
	\author[]{Marat Akhmet\thanks{Corresponding Author Tel.: +90 312 210 5355, Fax: +90 312 210 2972, E-mail: marat@metu.edu.tr} }
	\author[]{Ejaily Milad Alejaily}
	\affil[]{Department of Mathematics, Middle East Technical University, 06800 Ankara, Turkey}

	\date{}
	\maketitle
	
	We develop a new definition of fractals which can be considered as an abstraction of the fractals determined through self-similarity. The definition is formulated through imposing conditions which are governed the relation between the subsets of a metric space to build a porous self-similar structure. Examples are provided to confirm that the definition is satisfied by large class of self-similar fractals. The new concepts create new frontiers for fractals and chaos investigations.

	\section{Introduction}
	Fractals are class of complex geometric shapes with certain properties. One of the main features of the objects is self-similarity which can be defined as the property whereby parts hold similarity to the whole at any level of magnification \cite{Addison}. Fractional dimension is suggested by Mandelbrot to be a property of fractals when he defined a fractal as a set whose Hausdorff dimension strictly larger than its topological dimension \cite{Mandelbrot1}. Roots of the idea of self-similarity date back to the 17th century when Leibniz introduced the notions of recursive self-similarity \cite{Zmeskal}. The first mathematical definition of a self-similar shape was introduced in 1872 by Karl Weierstrass during his study of functions that were continuous but not differentiable. The most famous examples of fractals that display exact self-similarity are Cantor set, Koch curve and Sierpinski gasket and carpet which where discovered by Georg Cantor in 1883, Helge von Koch in 1904 and Waclaw Sierpinski in 1916 respectively. Julia sets, discovered by Gaston Julia and Pierre Fatou in 1917-19, gained significance in being generated using the dynamics of iterative function. In 1979, Mandelbrot visualized Julia sets including the most popular fractal called Mandelbrot set.
	
	In this paper we introduce a new mathematical concept and call it abstract fractal. This concept is an attempt to establish a pure foundation for fractals by abstracting the idea of self-similarity. We define the abstract fractal as a collection of points in a metric space. The points are represented through an iterative construction algorithm with specific conditions. The conditions are introduced to governor the relationship between the sets at each iteration. Our approach of construction is based on the concept of porosity rather than the roughness notion introduced by Mandelbrot. Porosity is an intrinsic property of materials and it is usually defined as the ratio of void volume to total volume \cite{Anovitz}. The concept of porosity plays an important role in several fields of research such as geology, soil mechanics, material science, civil engineering, etc. \cite{Anovitz,Ganji}. Fractal geometry has been widely used to study properties of porous materials. However, the concept of porosity was not utilized as as a criterion for fractal structures, and the relevant researches have investigated the relationship between porosity and fractalness \cite{Davis,Yu,Huang,Guyon,Cai}. For instance several researches such as \cite{Puzenko,Tang,Xia} determined fractal dimension of some pore-structures using the pore properties of them. The simplicity and importance of the porosity concept insistently invite us to develop a new definition of fractals through porosity. In other words the property should be involved in fractal theory as a feature to be equivalent to self-similarity and fractional dimension. This needs to specify the concept of porosity to surfaces and lines. In the present paper, we do not pay attention to equivalence between the definition of fractals in terms of porosity and those through self-similarity and dimension rather we introduce an abstract definition which, we hope, to be useful in application domains.  
		
	\section{The Definition} \label{AbsFractals}

	In this paper we shall consider the metric measure space defined by the triple $ (X, d, \mu) $, where $ (X, d) $ is a compact metric space, $ d $ is a metric on $ X $ and $ \mu $ is a measure on $ X $.

	To construct abstract fractal, let us consider the initial set $ F \subset X $ and fix two natural numbers $ m $ and $ M $ such that $ 1 < m < M $. We assume that there exist $ M $ nonempty disjoint subsets, $ F_{i}, \; i=1, 2, ... M $, such that $ \mathcal{F} = \cup_{i=1}^{M} \mathcal{F}_i $. For each $ i=1, 2, ... m $, again, there exist $ M $ nonempty disjoint subsets $ F_{ij}, \; j=1, 2, ... M $ such that $ \mathcal{F}_i = \cup_{j=1}^{M} \mathcal{F}_{ij} $. Generally, for each $ i_1, i_2, ..., i_n, \; i_k=1, 2, ... m $, there exist $ M $ nonempty disjoint sets $ F_{i_1 i_2 ... i_nj}, \; j=1, 2, ... M $, such that $ \mathcal{F}_{i_1 i_2 ... i_n} = \cup_{j=1}^{M} \mathcal{F}_{i_1 i_2 ... i_n j} $, for each natural number $ n $. The following conditions are needed: 
	
	\noindent There exist two positive numbers, $ r $ and $ R $, such that for each natural number $ n $ we have
	\begin{equation} \label{RatioCond}
	r \leq \frac{ \sum_{j=1}^{m} \mu \big( F_{i_1 i_2 ... i_{n-1} j} \big)}{ \sum_{j=m+1}^{M} \mu \big( F_{i_1 i_2 ... i_{n-1} j} \big)} \leq R.
	\end{equation}
	where $ i_k = 1, 2, ... , m, \; k=1, 2, ... n-1 $. We call the relation (\ref{RatioCond}) the \textit{ratio condition}. The numbers $ r $ and $ R $ in (\ref{RatioCond}) are characteristics for porosity. Another condition is the \textit{adjacent condition} and it is formulated as follows: 
	
	\noindent For each $ i_1 i_2 ... i_n, \; i_k=1, 2, ... , m $ there exists $ j, \; j=m+1, m+2, ... , M $, such that
	\begin{equation} \label{AdjtCond}
	d(F_{i_1 i_2 ... i_n}, F_{i_1 i_2 ... i_{n-1} j}) = 0.
	\end{equation}
	We call $ F_{i_1 i_2 ... i_{n-1} i_n} $ a complement set of order $ n $ if $ i_k=1, 2, ... , m, \; k=1, 2, ... n-1 $ and $ i_n = m+1, m+2, ... , M $. 
	
	An accumulation point of any couple of complement sets does not belong to any of them. We dub this stipulation the \textit{accumulation condition}.
		
	Let us define the diameter of a bounded subset $ A $ in $ X $ by $ \mathrm{diam}(A) = \sup \{ d(\textbf{x}, \textbf{y}) : \textbf{x}, \textbf{y} \in A \} $. Considering the above construction, we assume that the \textit{diameter condition} holds for the sets $ F_{i_1 i_2 ... i_n} $, i.e.,
	\begin{equation} \label{DiamCond}
	\max_{i_k=1,2, ... M} \mathrm{diam}(F_{i_1 i_2 ... i_n}) \to 	0 \;\; \text{as} \;\; n \to \infty.
	\end{equation}
	Fix an infinite sequence $ i_1 i_2 ... i_n ... \, $. The diameter conditions as well as the compactness of $ X $ imply that there exists a sequence $ (p_n) $, such that $ p_0 \in F $, $ p_1 \in F_{i_1} $, $ p_2 \in F_{i_1 i_2} $, ... , $ p_n \in F_{i_1 i_2 ... i_n}, \; n=1, 2, ...\, $, which converges to a point in $ X $. The points are denoted by $ F_{i_1 i_2 ... i_n ...} $.
	
	We define the \textit{abstract fractal} $ \mathcal{F} $ as the collection of the points $ F_{i_1 i_2 ... i_n ...} $ such that $ i_k=1,2, ... m $, that is
	\begin{equation} \label{AbsFracSet}
	\mathcal{F} =  \big\{F_{i_1 i_2 ... i_n ... } \; | \; i_k=1,2, ... m \big\}, 
	\end{equation}
	provided that the above four conditions hold. The subsets of $ \mathcal{F} $ can be represented by
	\begin{equation} \label{AbsFracSubSet}
	\mathcal{F}_{i_1 i_2 ... i_n} = \big\{ F_{i_1 i_2 ... i_n i_{n+1} i_{n+2} ... } \; | \; i_k=1, 2, ... , m \big\},
	\end{equation}
	where $ i_1 i_2 ... i_n $ are fixed numbers. We call such subsets \textit{subfractals} of order $ n $.
	
	\section{Abstract  structure of  geometrical fractals}
	
	In this section we find the pattern of abstract fractal in some geometrical well-known fractals, Sierpinski carpet, Pascal triangles and Koch curve.
		
	\subsection{The Sierpinski Carpet}
	
\begin{figure}[H]
	\centering
	\subfigure[]{\includegraphics[width = 1.6in]{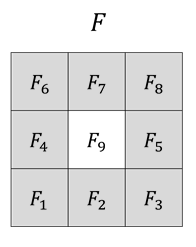}} \hspace{1cm}
	\subfigure[]{\includegraphics[width = 1.2in]{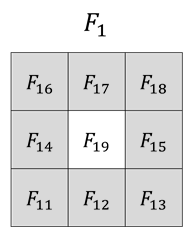}} \hspace{1cm} 
	\subfigure[]{\includegraphics[width = 1.5in]{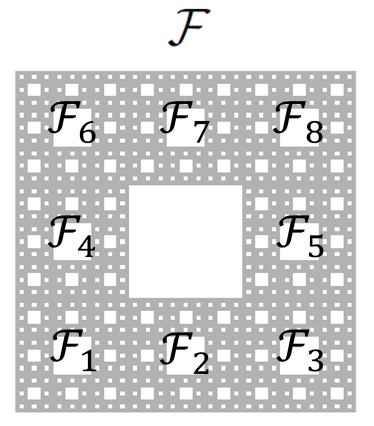}}
	\caption{Sierpinski carpet}
	\label{SCConst}   				
	\end{figure}

	To construct an abstract fractal corresponding to the Sierpinski carpet, let us consider a square as an initial set $ F $. Firstly, we divide $ F $ into nine $ (M=9) $ equal squares and denote them by $ F_{i}, \; i_1=1, 2, ... 9 $ (see Fig. \ref{SCConst} (a)). In the second step, each square $ F_{i}, \; i=1, 2, ... 8 $ is again divided into nine equal squares denoted as $ F_{i j}, \; j=1, 2, ... 9 $. Figure \ref{SCConst} (b) illustrates the sub-squares of $ F_1 $. We continue in this way such that at the $ n^{th} $ step, each set $ F_{i_1 i_2 ... i_{n-1}}, \; i_k=1, 2, ... 8 $, is divided into nine subset $ F_{i_1 i_2 ... i_{n-1}j}, \; j=1, 2, ... 9 $. For the Sierpinski carpet the number $ m $ is $ 8 $, and the measure ratio (\ref{RatioCond}) can be evaluated as follows. If we consider the first order sets $ F_{i_1}, \; i_1=1, 2, ... 9 $, then
	\[ \frac{ \sum_{j=1}^{8} \mu \big( F_{j} \big)}{ \mu \big( F_{9} \big)} = 8. \]
	
	Thus, the ratio condition holds. From the construction, we can see that each $ F_{i_1 i_2 ... i_n}, \; i_k=1, 2, ... , 8 $ has common boundary with $ F_{i_1 i_2 ... i_{n-1} j}, \; j= 9 $. Therefore, the adjacent condition holds. Since the construction consists of division into smaller parts, the diameter condition is also valid. Moreover, It is clear that the accumulation condition holds as well.
	
	As a result, the points of the desired abstract fractal $ \mathcal{F} $ can be represented as $ F_{i_1 i_2 ... i_n ... } $ and the abstract Sierpinski carpet is defines by
	\[ \mathcal{F} =  \big\{F_{i_1 i_2 ... i_n ... } | i_k=1,2, ... 5 \big\}. \]
	Figure \ref{SCConst} (c) shows the set $ \mathcal{F} $ and illustrates its $ 1^{st} $ order subfractals.
	
	\subsection{The Pascal Triangle}
	
\begin{figure}[H]
	\centering
	\subfigure[]{\includegraphics[width = 2.0in]{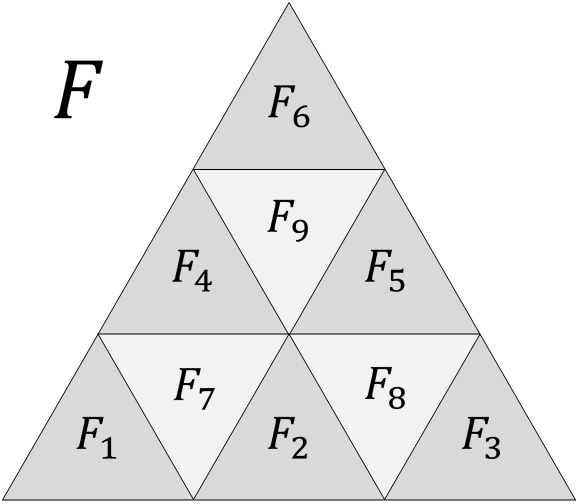}} \hspace{1cm}
	\subfigure[]{\includegraphics[width = 1.6in]{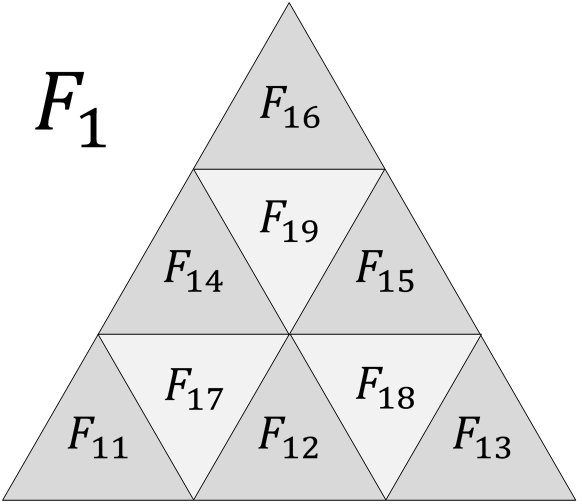}} \hspace{1cm}
	\subfigure[]{\includegraphics[width = 2.0in]{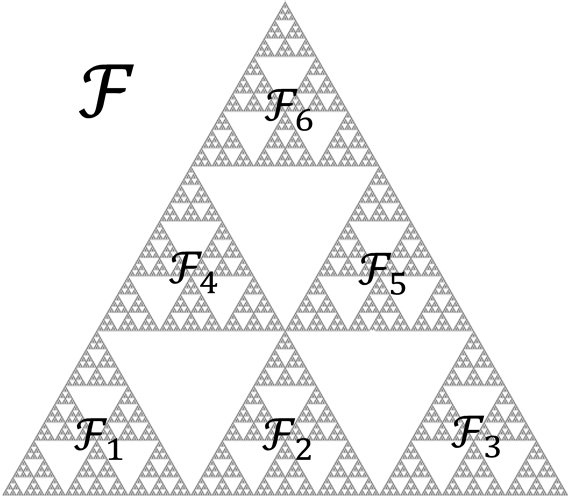}}	
	\caption{Pascal triangle modulo 3}
	\label{PTConst}   				
	\end{figure}

	Pascal triangle is a mathematical structure consists of triangular array of numbers. Triangular fractals can be obtained if these numbers are plotted using specific moduli. The Sierpinski gasket, for instance, is the Pascal's triangle modulo 2. Let us build an abstract fractal on the basis of a fractal associated with Pascal triangle modulo 3. Consider an equilateral triangle as an initial set $ F $. In the first step, we divide $ F $ into nine smaller equilateral triangles ad denote them as $ F_i, \; i=1, 2, ... , 9 $ as shown in Fig. \ref{PTConst} (a). Next, each triangle $ F_i, \; i=1, 2, ... , 6 $ is again into nine equilateral triangles named as $ F_{ij}, \; j=1, 2, ... , 9 $. Figure \ref{PTConst} (a) illustrates the second step for the set $ F_1 $. Similarly, the subsequent steps are performed such that at the $ n^{th} $ step, each set $ F_{i_1 i_2 ... i_{n-1}}, \; i_k=1, 2, ... 6 $, is divided into nine subset $ F_{i_1 i_2 ... i_{n-1}j}, \; j=1, 2, ... 9 $. In this case we have $ m=6 $ and $ M=9 $. Therefore,
	\[ \frac{ \sum_{j=1}^{6} \mu \big( F_{j} \big)}{ \sum_{j=7}^{9} \mu \big( F_{j} \big)} = 2, \]
	and the ratio condition holds. One can also verify that the adjacent, the accumulation, and the diameter conditions are also valid. Based on this, the points of the fractal can be defined by $ F_{i_1 i_2 ... i_n ... } $, and thus, the abstract Pascal triangle is defined by
	\[ \mathcal{F} =  \big\{F_{i_1 i_2 ... i_n ... } | i_k=1,2, ... 6 \big\}, \]
	and the $ n^{th} $ order subfractals can be written as
	\[ \mathcal{F}_{i_1 i_2 ... i_n} = \big\{ F_{i_1 i_2 ... i_n i_{n+1} ... } \; | \; i_k=1, 2, ... , 6 \big\}, \]
	where $ i_1 i_2 ... i_n $ are fixed numbers.

	\subsection{The Koch Curve}
	
\begin{figure}[H]
	\vspace{0.3cm}
	\centering
	\subfigure[]{\includegraphics[width = 2.6in]{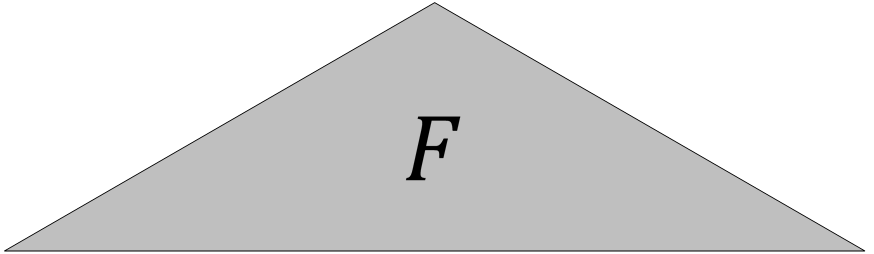}} \hspace{0.8cm}
	\subfigure[]{\includegraphics[width = 2.6in]{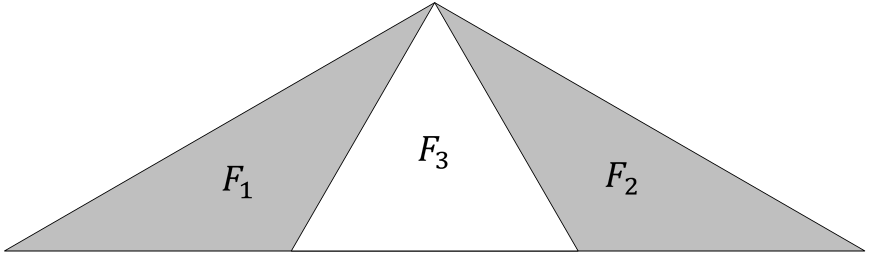}} \hspace{0.8cm}
	\subfigure[]{\includegraphics[width = 2.6in]{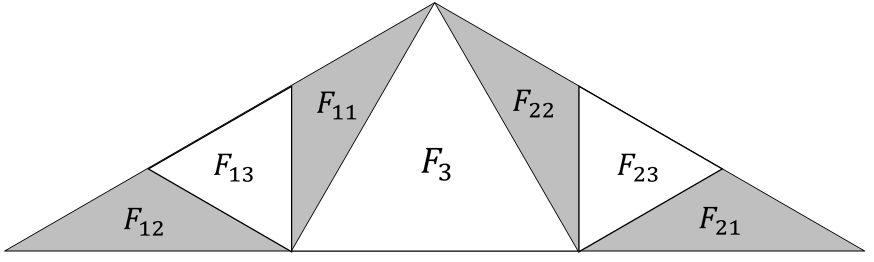}}
	\caption{Abstract Koch curve construction}
	\label{AKCConst}
	\vspace{0.3cm}  				
	\end{figure}

	In this subsection, we shall show how to build an abstract fractal $ \mathcal{F} $ conformable to the Koch curve. For this purpose, we consider the following construction of the Koch curve. Start with with an isosceles triangle $ F $ with base angles of $ 30^\circ $. The first step of the construction consists in dividing $ F $ into three equal-area triangles $ F_1, F_2 $ and $ F_3 $ (see Fig. \ref{AKCConst} (b)). The triangles $ F_1 $ and $ F_2 $ are isosceles with base angles of $ 30^\circ $, whereas the central triangle $ F_3 $ is an equilateral one. In the second step, each  $ F_i, \, i=1, 2 $ is similarly divided into three triangles, two isosceles, $ F_{i1} $ and $ F_{i2} $, and one equilateral, $ F_{i3} $. Figure \ref{AKCConst} (c) illustrate the step. In each subsequent step, the same procedure is repeated for each isosceles triangles resulting from the preceding step. That is, in the $ n^{th} $ step, each $ F_{i_1 i_2 ... i_{n-1}}, \; i_k=1, 2 $, is divided into three parts, two isosceles triangles $ F_{i_1 i_2 ... i_{n-1}j}, \; j=1, 2 $, with base angles of $ 30^\circ $, and one equilateral triangle $ F_{i_1 i_2 ... i_{n-1}3} $. In this construction, we have $ m=2 $ and $ M=3 $, thus, the measure ratio is
	\[ \frac{\mu (F_{1}) + \mu (F_{2})}{\mu (F_{3})} = 2, \]
	and the ratio condition holds. From the construction, it is clear that the adjacent, the accumulation, and  the diameter conditions are also valid. Based on this, the points in $ \mathcal{F} $ can be represented $ F_{i_1 i_2 ... i_n ...} $, and thus, the abstract Koch curve is defined by
	\begin{equation*} \label{AbsKochFract}
	\mathcal{F} =  \big\{F_{i_1 i_2 ... i_n ... } \; | \; i_k=1,2 \big\}. 
	\end{equation*}
	
	The $ n^{th} $ order subfractals of $ \mathcal{F} $ are represented by
	\begin{equation} \label{AbsKochSubFract}
	\mathcal{F}_{i_1 i_2 ... i_n} = \big\{ F_{i_1 i_2 ... i_n i_{n+1} i_{n+2} ... } \; | \; i_k=1, 2 \big\},
	\end{equation}
	where $ i_1 i_2 ... i_n $ are fixed numbers. Figure illustrats examples of $ 1^{st}, 2^{nd}, 3^{rd} $ and $ 4^{th} $ subfractals of the abstract Koch curve.

	\begin{figure}[H]
	\centering
	\includegraphics[width=0.50\linewidth]{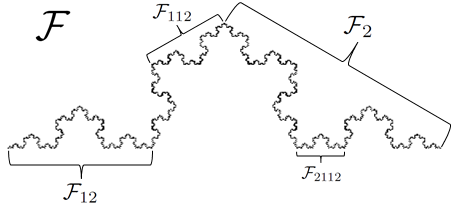}
	\caption{Subfractals of the abstract Koch curve}
	\label{KochSubFract}
	\end{figure}

	\section{Abstract self-similarity and Chaos}	

	In paper \cite{AkhmetSimilarity}, we have introduced the notion of the abstract self-similarity and defined a self-similar set by
	\begin{equation} \label{AbstSelfSimiSet}
	\mathcal{F} =  \big\{\mathcal{F}_{i_1 i_2 ... i_n ... } : i_k=1,2, ..., m, \; k=1, 2, ... \big\}, 
	\end{equation}
	where $ \mathcal{F}_{i_1 i_2 ... i_n ... }, \; i_k=1,2, ..., m $ represent the points of the set. For fixed indexes $ i_1, i_2, ..., i_n $, the subsets are expressed as
	\begin{equation} \label{AbstSelfSimiSubSet}
	\mathcal{F}_{i_1 i_2 ... i_n} = \bigcup_{j_k=1,2, ..., m } \mathcal{F}_{i_1 i_2 ... i_n j_1 j_2 ... },
	\end{equation}
	such that $ \mathcal{F}_{i_1 i_2 ... i_n} = \cup_{j=1}^{m} \mathcal{F}_{i_1 i_2 ... i_n j} $, for each natural number $ n $, where all sets $ \mathcal{F}_{i_1 i_2 ... i_n j}, \; j=1, 2, ..., m $, are nonempty, disjoint and satisfy the diameter condition.
	
	Based on the definition of the abstract self-similar set, we see that every abstract fractal is an abstract self-similar set, but the reverse is not necessarily valid.
	
	A similarity map $ \varphi $ for the abstract fractal $ \mathcal{F} $ can be defined by
	\[ \varphi(\mathcal{F}_{i_1 i_2 i_3 ...}) = \mathcal{F}_{i_2 i_3 i_4 ...}. \]
	
	Let us assume that the separation condition holds, that is, there exist a positive number $ \varepsilon_0 $ and a natural number $ n $ such that for arbitrary $ i_1 i_2 ... i_n $ one can find $ j_1 j_2 ... j_n $ so that
	\[ d \big( \mathcal{F}_{i_1 i_2 ... i_n} \, , \, \mathcal{F}_{j_1 j_2 ... j_n} \big) \geq \varepsilon_0, \]
	where $ \varepsilon_0 $ is the separation constant. Considering the results on chaos for self-similar set provided in \cite{AkhmetSimilarity}, it can be proven that the similarity map $ \varphi $ possesses the three ingredients of Devaney chaos, namely density of periodic points, transitivity and sensitivity. Moreover, $ \varphi $ possesses Poincar\`{e} chaos, which characterized by unpredictable point and unpredictable function \cite{AkhmetUnpredictable,AkhmetPoincare}. In addition to the Devaney and Poincar\`{e} chaos, it can be shown that the Li-Yorke chaos also takes place in the dynamics of the map. These results are summarized in the next theorem which can be proven in the similar way that explained in \cite{AkhmetSimilarity}. 
	
	\begin{theorem} \label{Thm1}
		If the separation condition holds, then the similarity map possesses chaos in the sense of Poincar\'{e}, Li-Yorke and Devaney..
	\end{theorem}
	
	That is the triple $ (\mathcal{F}, d, \varphi) $ is a self-similar space and $ \varphi $ is chaotic in the sense of Poincar\'{e}, Li-Yorke and Devaney.
		
	\section{Abstract Fractals and Iterated Function System}

	Iterated function system (IFS) is a powerful tool for the construction of fractal sets. It is defined by a family of contraction mappings $ w_n, n=1, 2, ... \, N $ on a complete metric space $ (X, d) $ \cite{Hutchinson,Barnsley}. The procedure starts with choosing an initial set $ A_0 \in \mathcal{B}(X) $, where $ \mathcal{B}(X) $ is the space of the non-empty compact subsets of $ X $, then iteratively applying the map $ W=\{ w_n, n=1, 2, ... \, N \} $ such that $ A_{k+1}= W(A_k) = \bigcup_{n=1}^N A_k^n $, where $ A_k^n = w_n(A_k) $. The fixed point of this map, $ A = W(A) = \lim_{k \to \infty} W^k(A_0) \in \mathcal{B}(X) $, is called the attractor of the IFS which represents the intended fractal.

\begin{figure}[H]
	\centering
	\subfigure[]{\includegraphics[width = 2.8in]{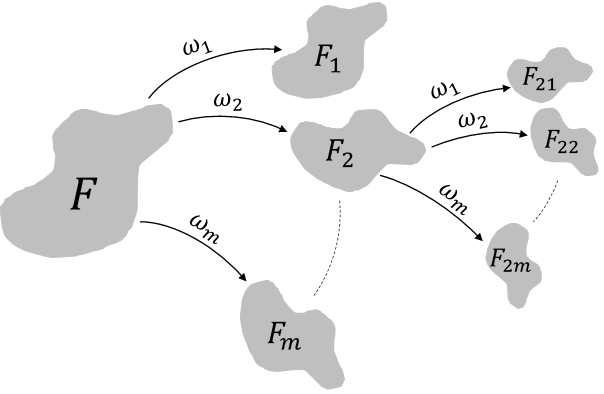}} \hspace{1.5cm}
	\subfigure[]{\includegraphics[width = 2.8in]{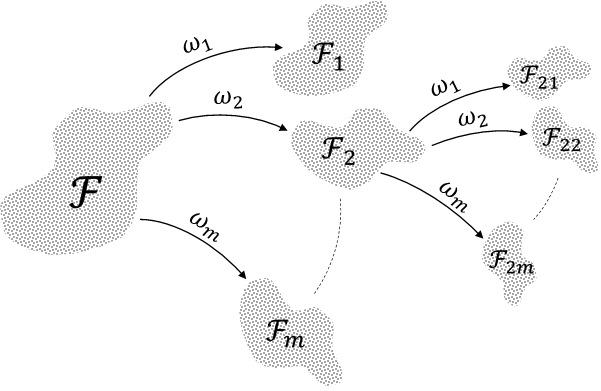}}
	\caption{IFS}
	\label{IFS}
	\end{figure}

	The idea of the structure of the abstract fractal can be realized using IFS. The fractal constructed by an IFS is an invariant set. Therefore, the subsets at each step of constructions can be determined using the maps $ w_n, n=1, 2, ... \, m $ as illustrated if Fig. \ref{IFS} (a). Similarly, the maps transform each subfractal into subfractals of the subsequent order. Figure \ref{IFS} (b) demonstrates the action of $ w_n $'s on the abstract fractal. The difference between this case and the above IFS fractal construction is that the sets $ \mathcal{F}_{i_1 i_2 ... i_n} $ are fractals in themselves, whereas the sets $ A_k^n $ are not.

	Utilizing the idea, moreover, each subfractal can be expressed in terms of the iterated images of whole fractal $ \mathcal{F} $, that is
	\[ \mathcal{F}_i= w_i (\mathcal{F}), \; \mathcal{F}_{ij}= w_j (w_i (\mathcal{F})), \; F_{ijk}= w_k (w_j (w_i (\mathcal{F}))), ... \, , \]
	thus, in general we have
	\[ \mathcal{F}_{i_1 i_2 ... i_n}= w_{i_n} (w_{i_{n-1}} ( ... (w_{i_1} (\mathcal{F}))) ... ), \]
	from which we can define a point belong to the fractal as the limit of the iterated images of $ \mathcal{F} $,
	\[ \mathcal{F}_{i_1 i_2 ...}= \lim_{n \to \infty} w_{i_n} (w_{i_{n-1}} ( ... (w_{i_1} (\mathcal{F}))) ... ). \]

	The existing of a separation constant $ \varepsilon_0 $ can be expressed in terms of $ w_n $ such that the condition is satisfied if
	\[ \min_n \inf_{\substack{i_n, j_n= \\ 1, 2 , ... \, m}} d \big( w_{i_n} ( ... (w_{i_1} (\mathcal{F})),  w_{j_n} ( ... (w_{j_1} (\mathcal{F})) \big) \geq \varepsilon_0. \]

	In addition to the construction of fractals, the IFS is used to prove chaos for the so-called totally disconnected IFS corresponding to certain classes of self-similar fractals like the Cantor set \cite{BarnsleyB}. The proof consists of construction of a dynamical system $ \{A; S\} $, where $ S: A \to A $ is the shift transformation defined by $ S(a) = w_n^{-1}(a) $ for $ a \in W_n(A) $. The system is called the shift dynamical system associated with the IFS, and then it showed to be topologically conjugate to the shift map on the $ N $-symbols code space. We see that this approach follows the usual construction of chaos which begins with defining a map with certain properties where the conjugacy to a well known chaotic map is the major key in discovering the chaotic nature of fractals. Again we emphasize that this approach is only applicable for the totally disconnected fractals namely the well-known Cantor sets. Differently, our approach is characterized by the similarity map $ \varphi: \mathcal{F} \to \mathcal{F} $ which, with regard to IFS approach, can be seen as an abstraction of the geometric essence of the transformation $ S $. Using the idea of indexing the domain elements allows to define the abstract map $ \varphi $. This has shortened the way of chaos proving by eliminating the need for topological conjugacy. Moreover, it becomes possible to investigate the chaotic nature in several classes of fractals such as the Sierpinski fractals and the Koch curve.
	
	\section{Discussion}
	
	The fractal concept is axiomatically linked with the notion self-similarity. This is why it is considered to be one of the two acceptable definitions of fractals. That is, a fractal can be defined as a set that display self-similarity at all scales. Mandelbrot define a fractal as a set whose Hausdorff dimension strictly larger than its topological dimension. In the present research, we introduce a conception of abstract fractal which can be considered as another criterion of fractalness. Indeed, the idea of the abstract fractal centers around the self-similarity property and many self-similar fractals like the Cantor sets and the Sierpinski fractals are shown to be fractals in the sense of the abstract fractal. These fractals are also satisfied the Mandelbrot definition. Moreover, in our previous paper \cite{AkhmetSimilarity}, we have also shown that the set of symbolic strings satisfies the definition of abstract self-similarity. Because of these facts, we believe that the notion of abstract fractal deserves to be the third definition of fractals and we hope it will be accepted by the mathematical community. Considering the abstract fractal as a new definition of fractal may open new opportunities for more theoretical investigations in this field as well as new possible applications in science and engineering. For example, we may start with the equivalency between these definitions. It is known that the fractals that display exact self-similarity at all scales satisfy Mandelbrot definition of fractal. The proposed definition satisfies the self-similarity since it is the main pivot of the concept of the abstract fractal. But the interesting question is: Does the abstract fractal agree with Mandelbrot definition? The notion of Hausdorff dimension for the abstract fractal is not yet developed enough to provide an answer to the question. However, the determination of the fractal dimension can possibly be performed based on two important properties. The first one is the self-similarity of the abstract fractal which may provide a self-similar dimension that can be assumed to be equivalent to the Hausdorff dimension. The second one is the accumulation condition combining perhaps with the diameter condition. These properties are essential for describing the geometry of fractals, therefore, the fractal dimension can be characterized in terms of them. This is why the definition in our paper can give opportunities to compare abstract fractals with fractals defined through dimension. Furthermore, the suggested fractal definition can be elaborated through chaotic dynamics development, topological spaces, physics, chemistry, neural network theories development.

\end{document}